\documentclass[10pt,a4paper,draft]{article}

\usepackage{bbm}
\usepackage{pifont}
\usepackage{bbding}
\usepackage[leqno]{amsmath}
\usepackage{amsmath, esint}
\usepackage{mathrsfs}
\usepackage{amsfonts}
\usepackage{amssymb}
\usepackage{amsfonts,amssymb,amsmath,indentfirst,cases,amsthm}
\usepackage{epsfig}
\usepackage[numbers,sort&compress]{natbib}
\usepackage[center]{titlesec}
\usepackage{fancyhdr}

\allowdisplaybreaks
\def\ess~sup{\mathop{\rm ess~sup}}

\numberwithin{equation}{section}

\newenvironment{key words}{\emph{\texttt{Keywords}}\mbox{  }}{ }
\newtheorem{theorem}{Theorem}[section]
\newtheorem{lemma}[theorem]{Lemma}
\newtheorem{remark}[theorem]{Remark}
\newtheorem{corollary}[theorem]{Corollary}
\newtheorem{proposition}[theorem]{Proposition}

\theoremstyle{remark}

\theoremstyle{plain}

\makeatletter

\newcommand{\Rmnum}[1]{\expandafter\@slowromancap\romannumeral #1@}
\makeatother
\begin{document}
\pagestyle{fancy}
\lhead{}
\chead{Junli Zhang, Pengcheng Niu}

\rhead{}

\renewcommand{\headrulewidth}{0.6pt}
 \title{\textbf{On a Conjecture of Cai-Zhang-Shen \\
 for Figurate Primes
}}
\author{Junli Zhang$^1$, Pengcheng Niu$^2$\thanks{Corresponding
author's E-mail: pengchengniu@nwpu.edu.cn(P. Niu)}\\
\small{1. School of Mathematics and Data Science, Shaanxi University of Science and Technology,}\\
\small{ Xi'an, Shaanxi, 710021, P. R. China}\\
\small{2. School of Mathematics and Statistics, Northwestern Polytechnical
University,}\\ \small{ Xi'an, Shaanxi, 710129, P. R. China}}
\date{}
\maketitle

\maketitle {\bf Abstract}\
A conjecture of Cai-Zhang-Shen for figurate primes says that every integer $k>1$ is the sum of two figurate primes. In this paper we give an equivalent proposition to the conjecture. By considering extreme value problems with constraints about the conjecture in the cases of odd and even integers and using the method of Lagrange multipliers, Cardano formula for cubic equations and the contradiction, we prove the conjecture.
\section{Introduction}

From 18th century it is known the so called Goldbach's binary conjecture which says that every even number greater than 2 can be written as the sum of two primes. This problem payed attention of many mathematicians, and, unfortunately, it is not solved till our days, see Apostol \cite{A76}, Chen \cite{C73}, Oliveira e Silva, Herzog and Pardi \cite{OHP14}, Pan and Pan \cite{PP92}, Wang \cite{W84}.

A binomial coefficient of the form $\left( {_s^{{p^r}}} \right)$ is called a figurate prime, where $p$ is a prime, $r \ge 1$ and $s \ge 0$ are integers. The collection of figurate primes includes 1, all primes and their powers, see \cite{C14}. It is well known that numbers of figurate primes and usual primes not larger than $x$ own the same density. In 2015, Cai, Zhang and Shen in \cite{CZS15} proposed a conjecture (we call it Cai-Zhang-Shen conjecture):
$$\hbox{every \;integer}\;k>1  \;\hbox{is\; the \;sum \;of \;two\; figurate \;primes},$$
and pointed out that the conjecture is true for integers up to ${10^7}$. In this paper we will discuss the conjecture and confirm that it is true.

Denote the characteristic function of figurate primes $i$ by $\delta (i)$, i.e., $\delta (i) = 1,\;\hbox{ when }\;i\;\hbox{is \; a\; figurate \;prime};$ $\delta (i) = 0,\;\hbox{ when }\;i\;\hbox{is\; not \;a\; figurate \;prime}.$
We claim that the Cai-Zhang-Shen conjecture for every integer $k \ge 3$ is equivalent to
\begin{equation}\label{eq11}
\sum\limits_{i = 1}^{k - 1} {\delta (i)\delta (k - i)}  > 0,\;\;k \ge 3.
\end{equation}
In fact, if \eqref{eq11} holds, then there exists $i$ such that
\[\delta (i)\delta (k - i) > 0,\]
that is $\delta (i) = \delta (k - i) = 1$, which implies that $i$ and $k-i$ are figurate primes, and the sum is $k$. Cai-Zhang-Shen conjecture is true. Conversely, if Cai-Zhang-Shen conjecture is true, that is every integer $k$ can be expressed as the sum of two figurate primes $i$ and $k-i$, then $\delta (i)\delta (k - i) > 0$ by $\delta (i) = \delta (k - i) = 1$, i.e., \eqref{eq11} is proved.

We can also give the equivalent descriptions for odd and even integers respectively. Let
$$I_n = \{i| \delta(i) = 1 \;\hbox{for\; some\; integer}\; i=1,2,\dots, {n - 1}\},$$
and by $l$ the number of figurate primes not being greater than $n-1$. We always let
\[l > {10^4}.\]
For odd integer $k = 2n - 1$, we take $N > 2n - 1$ satisfying $\delta (N) = 0$. Then Cai-Zhang-Shen conjecture is equivalent to
\begin{equation}\label{eq12}
\sum\limits_{i \in {I_n}} {\delta (i)\delta (2n - 1 - i)}  + \delta {(N)^2} > 0,\;n\ge 3.
\end{equation}
For even integer $k = 2n$, Cai-Zhang-Shen conjecture is equivalent to
\begin{equation}\label{eq13}
\sum\limits_{i \in {I_n}} {\delta (i)\delta (2n - i)}  + \delta {(n)^2} > 0,\;n\ge 3.
\end{equation}

The main result of the paper is
\begin{theorem}\label{Th1}
Cai-Zhang-Shen conjecture is ture.
\end{theorem}

We will divide odd integers and even integers to prove Theorem \ref{Th1}. The detailed  proof is given only in the case of odd integers, which can be similarly obtained in the case of even integers. Based on the properties satisfied by the characteristic function of the figurate primes, we introduce the objective function $f(x)$ ($x \in {\mathbb{R}^{2l + 1}}$), and two constraints $g(x) = 0$ and $h(x) = 0$. By testing that the set $A$ constructed by constraints is bounded, and the Jacobi determinant of two functions $g(x)$ and $h(x)$ is not 0, and then using the method of Lagrange multipliers, one shows $f(x) > 0$ on the set $A$. Under the assumption that Cai-Zhang-Shen conjecture is not true, the contradiction is obtained.

We emphasize the difficulties here, one is how to select the applicable objective function and constraints, especially the constraints, and the other is how to prove $f(x) > 0$ on $A$. Here the application of Cardano formula is successful.

Since Cai-Zhang-Shen conjecture is equivalent to \eqref{eq11}, we have from Theorem \ref{Th1} that
\begin{corollary}\label{Co1}
\eqref{eq11} holds.
\end{corollary}

This paper is organized as follows. The proof of Theorem \ref{Th1} (odd integers) is given in Section 2. We introduce the objective function $f(x)$ and two constraints $g(x)=0$ and $h(x)=0$. Using the method of Lagrange multipliers, one solves the minimum point of $f(x)$ on $A$ and infers $f(x)>0$ on $A$. Under the assumption that Cai-Zhang-Shen conjecture is not ture, the contradiction is derived. Therefore, Theorem \ref{Th1} (odd integers) is proved. Two propositions used in Section 2 are proved in Section 3. In Section 4, we prove theorem \ref{Th1} (even integers). Since the proof is similar to the previous sections, we only describe the related extreme value problem with constraints, and omit the details. Some conclusions are given in Section 5.

At the end of this section, let us state the method of Lagrange multipliers (e.g., refer to \cite{S12}) which will be used. For seeking the maximum and minimum values of $f(x)$($x \in {\mathbb{R}^n}$) with constraints
\[{g_i}(x) = 0\;(i = 1,2, \cdots ,m,m < n)\]
(assuming that these extreme values exist and the rank of Jacobian matrix
\[\frac{{\partial ({g_1}, \cdots ,{g_m})}}{{\partial ({x_1}, \cdots ,{x_n})}}\]
of ${g_i}(x)$ ($i = 1,2, \cdots ,m$) is $m$):\\
 (a) find all $x \in {\mathbb{R}^n},{\lambda _1}, \cdots ,{\lambda _m} \in \mathbb{R}$ such that
\[\frac{{\partial f}}{{\partial {x_i}}} + {\lambda _1}\frac{{\partial {g_1}}}{{\partial {x_i}}} +  \cdots  + {\lambda _m}\frac{{\partial {g_m}}}{{\partial {x_i}}} = 0,\;i = 1, \cdots ,n,\]
\[{g_i}(x) = 0,\;i = 1,2, \cdots ,m,\]
where $x$ is the stationary point and ${\lambda _1}, \cdots ,{\lambda _m}$ are multipliers;\\
(b) evaluate $f$ at all the points $x$ that result from (a). The largest of these values is the maximum value of $f$ and the smallest is the minimum value of $f$.
\section{\textbf{ Proof of Theorem \ref{Th1} (odd integers) }\label{Section 2}}

The following is Cardano formula for cubic equations:
\begin{lemma}\label{Le21}
Given the equation
$${y^3} + 3py + 2q = 0,$$
if $D = {p^3} + {q^2} > 0$, then there is a real solution
$$y = {u_ + } + {u_ - },$$
where
$${u_ + } = {\left( { - q + \sqrt D } \right)^{{\raise0.5ex\hbox{$\scriptstyle 1$}
\kern-0.1em/\kern-0.15em
\lower0.25ex\hbox{$\scriptstyle 3$}}}},\;{u_ - } = {\left( { - q - \sqrt D } \right)^{{\raise0.5ex\hbox{$\scriptstyle 1$}
\kern-0.1em/\kern-0.15em
\lower0.25ex\hbox{$\scriptstyle 3$}}}}.$$
\end{lemma}

\textbf{Proof of Theorem \ref{Th1} (odd integers)} Suppose that Cai-Zhang-Shen conjecture for odd integers is not true, namely there exists an odd integer $2n-1$ such that $2n-1$ can not be expressed as the sum of two figurate primes. Denote figurate primes not larger than $n-1$ by ${i_1},{i_2}, \cdots ,{i_l}\;({i_1} < {i_2} <  \cdots  < {i_l})$, and so
\[{I_n} = \left\{ {{i_1},{i_2}, \cdots ,{i_l}} \right\},\]
and $\delta ({i_1}) =1,\;  \delta ({i_2}) =1,\;\cdots  \; \delta ({i_l}) = 1;$ let
\[P = \left( {\delta ({i_1}), \cdots ,\delta ({i_l}),\delta (2n - 1 - {i_l}), \cdots ,\delta (2n - 1 - {i_1}),\delta (N)} \right),\]
i.e., components of $P$ are of \[\delta ({i_1}) =  \cdots  = \delta ({i_l}) = 1,\;\delta (2n - 1 - {i_1}) =  \cdots  = \delta (2n - 1 - {i_l}) = \delta (N) = 0,\]
Clrarly, $P \in {\mathbb{R}^{2l + 1}}$.

We introduce a function on ${\mathbb{R}^{2l+1}}$:
\begin{equation}\label{eq21}
f(x) = \sum\limits_{i \in I_n} {{x_i}{x_{2n - 1 - i}}}  + s{x_N^2},
\end{equation}
where \[s = \frac{3}{8}{l^{\frac{1}{3}}}.\]
Since $P$ satisfies
\[\sum\limits_{i \in {I_n}} {\left( {\delta {{(i)}^2}{\rm{ + }}\delta {{(2n - 1 - i)}^2}} \right)}  + \varepsilon \delta (N) = l,\]
\[\sum\limits_{i \in {I_n}} {\delta (i)\delta (2n - 1 - i)}  + \gamma \delta {(N)^3} + \frac{1}{2}\varepsilon \delta (N) = 0,\]
we define two functions on ${\mathbb{R}^{2l+1}}$:
\begin{equation}\label{eq22}
g(x) = \sum\limits_{i \in {I_n}} {\left( {x_i^2 + x_{2n - 1 - i}^2} \right)}  + \varepsilon {x_N} - l,
\end{equation}
\begin{equation}\label{eq23}
h(x) = \sum\limits_{i \in {I_n}} {{x_i}{x_{2n - 1 - i}}}  + \gamma x_N^3 + \frac{1}{2}\varepsilon {x_N},
\end{equation}
where
\[\varepsilon  = \frac{{3{l^{\frac{2}{3}}}}}{{4\sqrt 2 }},\gamma  =  - \frac{1}{{4\sqrt 2 }}.\]
Consider the extreme values of $f(x)$ with constraints
\begin{equation}\label{eq24}
g(x) = 0\;\hbox{and}\;h(x) = 0.
\end{equation}
Denote
\begin{equation}\label{eq25}
A = \left\{ {x \in {\mathbb{R}^{2l + 1}}|g(x) = 0,h(x) = 0} \right\}.
\end{equation}

We describe two propositions whose proofs will put in Section 3.
\begin{proposition}\label{Pro21}
 The set $A$ is bounded and closed in ${\mathbb{R}^{2l + 1}}$.
\end{proposition}

\begin{proposition}\label{Pro22}
The rank of the Jacobian matrix for functions $g(x)$ and $h(x)$ on $A$ is 2.
\end{proposition}

\begin{remark}\label{Re21}
 Under the assumption that CZS conjecture is not true, we see that $P \in {\mathbb{R}^{2l + 1}}$ belongs to $A$, because $P$ satisfies \eqref{eq24}.
\end{remark}

\begin{remark}\label{Re22}
By Proposition \ref{Pro22}, there are infinite points in $A$, since there are $2l-1$ independent variables in $A$.
\end{remark}

\begin{remark}\label{Re23}
If Cai-Zhang-Shen conjecture is not true, then
\[\delta (2n - 1 - i) = \delta (N) = 0\;(i \in {I_n})\]
and
\begin{equation}\label{eq26}
f(P) = \sum\limits_{i \in {I_n}} {\delta (i)\delta (2n - 1 - i) + } \delta {(N)^2} = 0.
\end{equation}
\end{remark}

We write the Lagrange function
\begin{equation}\label{eq27}
Q(x,\lambda ,\mu ) = f(x) + \lambda g(x) + \mu h(x),
\end{equation}
and use the method of Lagrange multipliers to find all stationary points of $f(x)$ on $A$, and then prove \begin{center}
$f(x) > 0$ at these points,
\end{center}
which show
\begin{center}
$f(x) > 0$ on $A$.
\end{center}

1) For $i \in I_n$, we have
\begin{equation}\label{eq28}
\left\{ {\begin{array}{*{20}{c}}
{{Q_{{x_i}}} = {x_{2n - 1 - i}} + 2\lambda {x_i} + \mu {x_{2n - 1 - i}} = 0,}\\
{{Q_{{x_{2n - 1 - i}}}} = {x_i} + 2\lambda {x_{2n - 1 - i}} + \mu {x_i} = 0,}
\end{array}} \right.
\end{equation}
i.e.,
\[\left\{ {\begin{array}{*{20}{c}}
{2\lambda {x_i} + (1 + \mu ){x_{2n - 1 - i}} = 0,}\\
{(1 + \mu ){x_i} + 2\lambda {x_{2n - 1 - i}} = 0.}
\end{array}} \right.\]
The determinant of coefficients is
\begin{equation}\label{eq29}
\left| {\begin{array}{*{20}{c}}
{2\lambda }&{1 + \mu }\\
{1 + \mu }&{2\lambda }
\end{array}} \right| = {\left( {2\lambda } \right)^2} - {\left( {1 + \mu } \right)^2},
\end{equation}
hence\\
$a_1$)  ${\left( {2\lambda } \right)^2} - {\left( {1 + \mu } \right)^2} \ne 0$, ${x_i} = {x_{2n - 1 - i}} = 0$;\\
for
\[{\left( {2\lambda } \right)^2} - {\left( {1 + \mu } \right)^2} = 0,\]
we have\\
$a_2$) $\lambda  \ne 0,$ $2\lambda  =  - (1 + \mu )$, ${x_i} - {x_{2n - 1 - i}} = 0$;\\
$a_3$) $\lambda  \ne 0,$ $2\lambda  =  1 + \mu $, ${x_i} + {x_{2n - 1 - i}} = 0$;\\
$a_4$) $\lambda  = 0,$ $\mu  =  - 1$, ${x_i}$ and ${x_{2n - 1 - i}} $ are arbitrary.

2) For $i = N$, we have ${Q_{{x_N}}} = 2s{x_N} + \varepsilon \lambda  + 3\gamma \mu x_N^2 + \frac{1}{2}\varepsilon \mu  = 0,$ so
\begin{equation}\label{eq210}
3\gamma \mu x_N^2 + 2s{x_N} + \varepsilon \lambda  + \frac{1}{2}\varepsilon \mu  = 0,
\end{equation}
its discriminant is
\begin{equation}\label{eq211}
\Delta  = {\left( {2s} \right)^2} - 12\gamma \mu \left( {\varepsilon \lambda  + \frac{1}{2}\varepsilon \mu } \right),
\end{equation}
therefore\\
$b_1$) $\mu  = 0,$ $2s{x_N} + \varepsilon \lambda  = 0$ and ${x_N} = \frac{{ - \varepsilon \lambda }}{{2s}}$;\\
$b_2$) $\mu  \ne 0,\;\Delta  = 0$, ${x_N} =  - \frac{s}{{3\gamma \mu }}$;\\
$b_3$) $\mu  \ne 0,\;\Delta  > 0$, ${x_N} = \frac{{ - 2s + \sqrt \Delta  }}{{6\gamma \mu }}$;\\
$b_4$) $\mu  \ne 0,\;\Delta  < 0$, ${x_N} = \frac{{ - 2s - \sqrt \Delta  }}{{6\gamma \mu }}$.
\begin{remark}\label{Re24}
Note that $P$ is not a stationary point. In fact, components of $P$ do not satisfy ${a_1}),{a_2}),{a_3})$. If $P$ satisfies ${a_4})$, it knows $\mu  =  - 1$, which contradicts to $\mu  =  0$ by $b_1$); it gives ${x_N} \ne 0$ by $b_2$), which contradicts to the component $\delta (N) = 0$ of $P$; if $P$ satisfies $b_3$), then ${x_N} = \frac{{ - 2s + \sqrt \Delta  }}{{6\gamma \mu }} = 0$ and $2s = \sqrt \Delta  ,$ so $ - 12\gamma  \mu \left( {\varepsilon \lambda  + \frac{1}{2}\varepsilon \mu } \right) = 0$ from \eqref{eq211}, but $ - 12\gamma  \mu \left( {\varepsilon \lambda  + \frac{1}{2}\varepsilon \mu } \right) =  - 6\gamma \varepsilon  \ne 0$ by $\lambda  = 0$ and $\mu  =  - 1$ in ${a_4})$, a contradiction; if $P$ satisfies $b_4$), then ${x_N} = \frac{{ - 2s - \sqrt \Delta  }}{{6\gamma \mu }} = 0$ and $2s =  - \sqrt \Delta  ,$ and ${(2s)^2} = \Delta ,$ it gets $ - 12\gamma \mu \left( {\varepsilon \lambda  + \frac{1}{2}\varepsilon \mu } \right) = 0$ by \eqref{eq211}, but $ - 12\gamma  \mu \left( {\varepsilon \lambda  + \frac{1}{2}\varepsilon \mu } \right) =  - 6\gamma \varepsilon  \ne 0$ by $\lambda  = 0$ and $\mu  =  - 1$ in ${a_4})$, a contradiction. Hence $P$ does not satisfy ${b_1}) - {b_4})$ which shows that $P$ is not a stationary point.
\end{remark}

Let us discuss all combinations of $a_1)$-$a_4)$ and $b_1)$-$b_4)$ and prove $f(x) > 0$ at all stationary points.

\textbf{Case ${a_1}),{b_1})$:} Note that ${x_i} = {x_{2n - 1 - i}} = 0$ ($i \in I_n$) from $a_1)$. Using
\[0 = g(x) = \varepsilon {x_N} - l,\]
it solves
\begin{equation}\label{eq212}
{x_N} = \frac{l}{\varepsilon } = \frac{{4\sqrt 2 }}{3}{l^{\frac{1}{3}}}.
\end{equation}
Since
\[0 = h(x) = \gamma x_N^3 + \frac{\varepsilon }{2}{x_N} = {x_N}\left( {\gamma x_N^2 + \frac{\varepsilon }{2}} \right),\]
we have
\begin{equation}\label{eq213}
{x_N} = 0\;\hbox{ or }\;x_N^2 =  - \frac{\varepsilon }{{2\gamma }} = \frac{3}{2}{l^{\frac{2}{3}}}.
\end{equation}
It is different from ${x_N}$ in \eqref{eq212}, a contradiction.

\textbf{Case ${a_1}),{b_2})$:} It leads to a contradiction as in Case ${a_1}),{b_1})$.

\textbf{Case ${a_1}),{b_3})$:} It leads to a contradiction as in Case ${a_1}),{b_1})$.

\textbf{Case ${a_1}),{b_4})$:} It leads to a contradiction as in Case ${a_1}),{b_1})$.

\textbf{Case ${a_2}),{b_1})$:}  Noting $2\lambda  =  - (1 + \mu )$ and ${x_i} = {x_{2n - 1 - i}}$ by ${a_2})$, and $\mu  = 0$ by ${b_1})$, we obtain $\lambda  =  - \frac{1}{2}$ and also by ${b_1})$,
\begin{equation}\label{eq214}
{x_N} = \frac{{ - \varepsilon \lambda }}{{2s}} = \frac{\varepsilon }{{4s}} = \frac{{{l^{\frac{1}{3}}}}}{{2\sqrt 2 }}.
\end{equation}
Applying ${x_i} = {x_{2n - 1 - i}}$, we see
\[0 = g(x) = 2\sum\limits_{i \in {I_n}} {x_i^2}  + \varepsilon {x_N} - l,\]
\[0 = 2h(x) = 2\sum\limits_{i \in {I_n}} {x_i^2}  + 2\gamma x_N^3 + \varepsilon {x_N},\]
and so
\[2\gamma x_N^3 + l = 0,\]
then
\begin{equation}\label{eq215}
{x_N} = {\left( {\frac{{ - l}}{{2\gamma }}} \right)^{\frac{1}{3}}} = \sqrt 2 {l^{\frac{1}{3}}}.
\end{equation}
It is different from ${x_N}$ in \eqref{eq214}, a contradiction.

\textbf{Case ${a_2}),{b_2})$:} In virtue of ${x_i} = {x_{2n - 1 - i}}$ by ${a_2})$, similarly to Case ${a_2}),{b_1})$, we have
\[{x_N} = {\left( {\frac{{ - l}}{{2\gamma }}} \right)^{\frac{1}{3}}} = \sqrt 2 {l^{\frac{1}{3}}}.\]
It follows
$$f(x) = \sum\limits_{i \in {I_n}} {x_i^2}  + sx_N^2 \ge sx_N^2 = \frac{3}{8}{l^{\frac{1}{3}}}{\left( {\sqrt 2 {l^{\frac{1}{3}}}} \right)^2} = \frac{3}{4}l > 0.$$

\textbf{Case ${a_2}),{b_3})$:} We use ${x_i} = {x_{2n - 1 - i}}$ to derive $f(x) > 0$ as in Case ${a_2}),{b_2})$.

\textbf{Case ${a_2}),{b_4})$:} We use ${x_i} = {x_{2n - 1 - i}}$ to derive $f(x) > 0$ as in Case ${a_2}),{b_2})$.

\textbf{Case ${a_3}),{b_1})$:} It gives $2\lambda  = 1 + \mu $ and ${x_i} =  - {x_{2n - 1 - i}}$ by ${a_3})$ and $\mu  = 0$ by ${b_1})$, then $\lambda  = \frac{1}{2}$ and by ${b_1})$,
\begin{equation}\label{eq216}
{x_N} = \frac{{ - \varepsilon \lambda }}{{2s}} = \frac{{ - \varepsilon }}{{4s}} =  - \frac{{{l^{\frac{1}{3}}}}}{{2\sqrt 2 }}.
\end{equation}
On the other hand, using ${x_i} =  - {x_{2n - 1 - i}}$, it yields
\[0 = g(x) = 2\sum\limits_{i \in {I_n}} {x_i^2}  + \varepsilon {x_N} - l,\]
\[0 = 2h(x) =  - 2\sum\limits_{i \in {I_n}} {x_i^2}  + 2\gamma x_N^3 + \varepsilon {x_N},\]
so
\begin{equation}\label{eq217}
x_N^3 + \frac{\varepsilon }{\gamma }{x_N} - \frac{l}{{2\gamma }} = 0.
\end{equation}
Since
\[3p = \frac{\varepsilon }{\gamma },\quad 2q =  - \frac{l}{{2\gamma }},\]
and
\[p = \frac{\varepsilon }{{3\gamma }},\quad q = \frac{{ - l}}{{4\gamma }},\]
we have from Lemma \ref{Le21} and
\[D = {p^3} + {q^2} = {\left( {\frac{\varepsilon }{{3\gamma }}} \right)^3} + {\left( {\frac{{ - l}}{{4\gamma }}} \right)^2} =  - {l^2} + 2{l^2} = {l^2},\]
\[{u_ + } = {\left( { - q + \sqrt D } \right)^{\frac{1}{3}}} = {\left( {\frac{l}{{4\gamma }} + l} \right)^{\frac{1}{3}}} = {\left( { - \sqrt 2 l + l} \right)^{\frac{1}{3}}},\]
\[{u_ - } = {\left( { - q - \sqrt D } \right)^{\frac{1}{3}}} = {\left( {\frac{l}{{4\gamma }} - l} \right)^{\frac{1}{3}}} = {\left( { - \sqrt 2 l - l} \right)^{\frac{1}{3}}}\]
that a real solution to \eqref{eq217} is
\begin{equation}\label{eq218}
{x_N} = {u_ + } + {u_ - } = \left( { - {{\left( {\sqrt 2  - 1} \right)}^{\frac{1}{3}}} - {{\left( {\sqrt 2  + 1} \right)}^{\frac{1}{3}}}} \right){l^{\frac{1}{3}}} \approx  - 2.087{l^{\frac{1}{3}}}.
\end{equation}
It is different from ${x_N}$ in \eqref{eq216}), a contradiction.

\textbf{Case ${a_3}),{b_2})$:} Noting ${x_i} + {x_{2n - 1 - i}} = 0$ by ${a_3})$, it follows as in Case ${a_3}),{b_1})$ that
\[{x_N} =  \left( { - {{\left( {\sqrt 2  - 1} \right)}^{\frac{1}{3}}} - {{\left( {\sqrt 2  + 1} \right)}^{\frac{1}{3}}}} \right){l^{\frac{1}{3}}} \approx  - 2.087{l^{\frac{1}{3}}}.\]
Using
\[0 = h(x) = \sum\limits_{i \in {I_n}} {{x_i}{x_{2n - 1 - i}}}  + \gamma x_N^3 + \frac{1}{2}\varepsilon {x_N},\]
it implies
\begin{align*}
  f(x) & =  - \gamma x_N^3 - \frac{1}{2}\varepsilon {x_N} + sx_N^2 \\
   &  =  - \left( {\frac{{ - 1}}{{4\sqrt 2 }}} \right){\left( { - 2.087} \right)^3}l - \frac{1}{2}\frac{{3{l^{\frac{2}{3}}}}}{{4\sqrt 2 }}\left( { - 2.087} \right){l^{\frac{1}{3}}} + \frac{3}{8}{l^{\frac{1}{3}}}{\left( { - 2.087} \right)^2}{l^{\frac{2}{3}}}\\
   & =  - \frac{{{{\left( {2.087} \right)}^3}}}{{4\sqrt 2 }}l + \frac{{3 \cdot \left( {2.087} \right)}}{{8\sqrt 2 }}l + \frac{{3 \cdot {{\left( {2.087} \right)}^2}}}{8}l = 2.087l\left( { - \frac{{{{\left( {2.087} \right)}^2}}}{{4\sqrt 2 }} + \frac{3}{{8\sqrt 2 }} + \frac{{3 \cdot \left( {2.087} \right)}}{8}} \right)\\
   & = 2.087l\left( { - 0.769 + 0.265 + 0.75} \right) = 2.087l\left( { - 0.769 + 1.015} \right) > 0.
\end{align*}

\textbf{Case ${a_3}),{b_3})$:} It follows $f(x) > 0$ as in Case ${a_3}),{b_2})$.

\textbf{Case ${a_3}),{b_4})$:} It follows also $f(x) > 0$ as in Case ${a_3}),{b_2})$.

\textbf{Case ${a_4}),{b_1})$:} It knows $\mu  =  - 1$ by ${a_4})$, which contradicts to $\mu  =  0$ by ${b_1})$.

\textbf{Case ${a_4}),{b_2})$:} Note by ${a_4})$ and ${b_2})$, we have $\mu  =  - 1$ and
\[{x_N} =  - \frac{s}{{3\gamma \mu }} = \frac{s}{{3\gamma }} =  - \frac{1}{{\sqrt 2 }}{l^{\frac{1}{3}}}.\]
Using
\[h(x) = 0,\]
it derives
\begin{align*}
   f(x)& =  - \gamma x_N^3 - \frac{1}{2}\varepsilon {x_N} + sx_N^2 =  - \frac{{ - 1}}{{4\sqrt 2 }}{\left( {\frac{{ - {l^{\frac{1}{3}}}}}{{\sqrt 2 }}} \right)^3} - \frac{1}{2}\frac{{3{l^{\frac{2}{3}}}}}{{4\sqrt 2 }}\left( {\frac{{ - {l^{\frac{1}{3}}}}}{{\sqrt 2 }}} \right) + \frac{3}{8}{l^{\frac{1}{3}}}{\left( {\frac{{ - {l^{\frac{1}{3}}}}}{{\sqrt 2 }}} \right)^2} \\
  &  =  - \frac{1}{{4\sqrt 2 }}\cdot \frac{l}{{2\sqrt 2 }} + \frac{{3l}}{{16}} + \frac{3}{8}\cdot \frac{l}{2} =  - \frac{l}{{16}} + \frac{{3l}}{{16}} + \frac{{3l}}{{16}} > 0.
\end{align*}

\textbf{Case ${a_4}),{b_3})$:} Notes $\lambda  = 0$ and $\mu  =  - 1$ by ${a_4})$ and so $\Delta  = 0$, which contradicts to $\Delta  > 0$ by ${b_3})$.

\textbf{Case ${a_4}),{b_4})$:} As in Case ${a_4}),{b_3})$, it also follows a contradiction.

Noting that $A$ is a bounded closed set in ${\mathbb{R}^{2l + 1}}$ and $f(x)$ is continuous in ${\mathbb{R}^{2l + 1}}$, we know that $f(x)$ achieves the minimum value on $A$. Summing up above discussions, we  indeed prove that the minimum of $f(x)$ on $A$ is positive, and so
\[f(x) > 0,\;x \in A.\]

\textbf{ End of Proof of Theorem \ref{Th1} (odd integers)} Since one supposes that Cai-Zhang-Shen conjecture is not true, it follows $f(x) > 0$ ($x \in A$) from the above analysis and so
\[f(P) > 0\]
because of $P \in A$. But it contradicts to \eqref{eq26}. Theorem \ref{Th1} (odd integers) is proved.

\section{\textbf{Proofs of Propositions \ref{Pro21} and \ref{Pro21}}\label{Section 3}}

\textbf{Proof of Proposition \ref{Pro21}} The closeness of $A$ in \eqref{eq25} is evident. We divide
two steps to prove that $A$ is bounded, i.e., first prove that when the set $\{ {x_N}\} $ constructed by components ${x_N}$ of $x \in A$ is bounded, it concludes that $A$ is bounded; next prove that the set $\{ {x_N}\} $ must be bounded by the contradiction.

\textbf{Step 1} Suppose that the set $\{ {x_N}\} $ is bounded, then there exists a constant
$C > 0,$ such that $\left| {{x_N}} \right| \le C.$ It uses $g(x) = 0$ to show
\[\sum\limits_{j = 1}^l {\left( {x_{{i_j}}^2 + x_{2n - 1 - {i_j}}^2} \right)}  + x_N^2 = x_N^2 - \varepsilon {x_N} + l \le {C^2} + \varepsilon C + l.\]
Hence $A$ is bounded.

\textbf{Step 2} Let us prove the boundedness of $\{ {x_N}\} $ by the contradiction. Assume that $\{ {x_N}\} $ is unbounded, then for any positive integer $\alpha ,$ there exists ${x_{N,\alpha }}$ in $\{ {x_N}\} $, such that $\left| {{x_{N,\alpha }}} \right| \ge \alpha .$ So ${x_{N,\alpha }} \to \infty $  as $\alpha \to \infty$. For convenience, we simply denote ${x_N} \to \infty $. It follows from $g(x) = 0$ that
\begin{equation}\label{eq31}
 - \varepsilon {x_N} + l = \sum\limits_{j = 1}^l {\left( {x_{{i_j}}^2 + x_{2n - 1 - {i_j}}^2} \right)}
\end{equation}
and ${x_N} \to  \infty$ should be
\[{x_N} \to  - \infty ,\]
so there exists one or several components in ${x_{{i_j}}},{x_{2n - 1 - {i_j}}}(j = 1,2, \cdots ,l)$ tending to $\infty$. We consider the following subcases.

1) If ${x_{{i_1}}} \to \infty $ and ${x_{2n - 1 - {i_1}}},{x_{{i_j}}},{x_{2n - 1 - {i_j}}}(j = 2, \cdots ,l)$ are bounded, then we have
$x_{{i_1}}^2 \to  + \infty $ and from \eqref{eq31} that
\begin{equation}\label{eq32}
x_{{i_1}}^2 =  - x_{2n - 1 - {i_1}}^2 - \sum\limits_{j = 2}^l {\left( {x_{{i_j}}^2 + x_{2n - 1 - {i_j}}^2} \right)}  - \varepsilon {x_N} + l: =  - \varepsilon {x_N} + {C_1},
\end{equation}
where ${C_1}$ is finite, so
\[\frac{{x_{{i_1}}^2}}{{ - {x_N}}} \to \varepsilon ,\quad \frac{{\left| {{x_{{i_1}}}} \right|}}{{{{\left( { - {x_N}} \right)}^{\frac{1}{2}}}}} \to \sqrt \varepsilon  .\]
It yields from $h(x) = 0$ that
\begin{equation}\label{eq33}
 {x_{{i_1}}}{x_{2n - 1 - {i_1}}} =  - \sum\limits_{j = 2}^l {{x_{{i_j}}}{x_{2n - 1 - {i_j}}}}  - \gamma x_N^3 - \frac{1}{2}\varepsilon {x_N}: =  - \gamma x_N^3 - \frac{1}{2}\varepsilon {x_N} + {C_2},
\end{equation}
where ${C_2}$ is finite.

When ${x_{2n - 1 - {i_1}}} = 0,$ we have by \eqref{eq33} that
\[0 =  - \gamma x_N^3 - \frac{1}{2}\varepsilon {x_N} + {C_2},\]
and the right hand side tends to $ - \infty $ (noting $\gamma  < 0$), a contradiction.

When ${x_{2n - 1 - {i_1}}} \ne 0,$ it follows from \eqref{eq33} to see
\begin{align*}
   \frac{{{x_{{i_1}}}{x_{2n - 1 - {i_1}}}}}{{{{\left( { - {x_N}} \right)}^{\frac{1}{2}}}}}& =  - \gamma \frac{{x_N^3}}{{{{\left( { - {x_N}} \right)}^{\frac{1}{2}}}}} - \frac{1}{2}\varepsilon \frac{{{x_N}}}{{{{\left( { - {x_N}} \right)}^{\frac{1}{2}}}}} + \frac{{{C_2}}}{{{{\left( { - {x_N}} \right)}^{\frac{1}{2}}}}} \\
   & = \gamma \frac{{{{\left( { - {x_N}} \right)}^3}}}{{{{\left( { - {x_N}} \right)}^{\frac{1}{2}}}}} + \frac{1}{2}\varepsilon \frac{{ - {x_N}}}{{{{\left( { - {x_N}} \right)}^{\frac{1}{2}}}}} + \frac{{{C_2}}}{{{{\left( { - {x_N}} \right)}^{\frac{1}{2}}}}} \to  - \infty ,\quad (\gamma  < 0)
\end{align*}
but the left hand side tends to $ \pm \sqrt \varepsilon  {x_{2n - 1 - {i_1}}},$ a contradiction.

2) If ${x_{{i_1}}} \to \infty $ and ${x_{2n - 1 - {i_1}}} \to \infty $ and ${x_{{i_j}}},{x_{2n - 1 - {i_j}}}(j = 2, \cdots ,l)$ are bounded, then
\[x_{{i_1}}^2 + x_{2n - 1 - {i_1}}^2 \to  + \infty ,\]
it shows by \eqref{eq31} that
\begin{equation}\label{eq34}
x_{{i_1}}^2 + x_{2n - 1 - {i_1}}^2 =  - \sum\limits_{j = 2}^l {\left( {x_{{i_j}}^2 + x_{2n - 1 - {i_j}}^2} \right)}  - \varepsilon {x_N} + l: =  - \varepsilon {x_N} + {C_3},
\end{equation}
where ${C_3}$ is finite, so
\[\frac{{x_{{i_1}}^2 + x_{2n - 1 - {i_1}}^2}}{{ - {x_N}}} \to \varepsilon  > 0.\]
It gives from $h(x) = 0$ that
\begin{equation}\label{eq35}
{x_{{i_1}}}{x_{2n - 1 - {i_1}}} =  - \sum\limits_{j = 2}^l {{x_{{i_j}}}{x_{2n - 1 - {i_j}}} - \gamma x_N^3 - \frac{1}{2}\varepsilon {x_N}: = }  - \gamma x_N^3 - \frac{1}{2}\varepsilon {x_N} + {C_4},
\end{equation}
where ${C_4}$ is finite. We have by \eqref{eq35} that
\[\frac{{{x_{{i_1}}}{x_{2n - 1 - {i_1}}}}}{{ - {x_N}}} =  - \gamma \frac{{x_N^3}}{{ - {x_N}}} - \frac{1}{2}\varepsilon \frac{{{x_N}}}{{ - {x_N}}} + \frac{{{C_4}}}{{ - {x_N}}} = \gamma \frac{{ - x_N^3}}{{ - {x_N}}} + \frac{1}{2}\varepsilon \frac{{ - {x_N}}}{{ - {x_N}}} + \frac{{{C_4}}}{{ - {x_N}}} \to  - \infty ,\]
then
\[ + \infty  \leftarrow \frac{{2\left| {{x_{{i_1}}}{x_{2n - 1 - {i_1}}}} \right|}}{{ - {x_N}}} \le \frac{{x_{{i_1}}^2 + x_{2n - 1 - {i_1}}^2}}{{ - {x_N}}} \to \varepsilon ,\]
a contradiction.

3) If ${x_{{i_1}}} \to \infty $ and ${x_{{i_2}}} \to \infty $ and ${x_{2n - 1 - {i_1}}},{x_{2n - 1 - {i_2}}},{x_{{i_j}}},{x_{2n - 1 - {i_j}}}(j = 3, \cdots ,l)$ are bounded, then
\[x_{{i_1}}^2 + x_{{i_2}}^2 \to \infty \]
and from \eqref{eq31},
\begin{equation}\label{eq36}
 x_{{i_1}}^2 + x_{{i_2}}^2 =  - \varepsilon {x_N} - x_{2n - 1 - {i_1}}^2 - x_{2n - 1 - {i_2}}^2 - \sum\limits_{j = 3}^l {\left( {x_{{i_j}}^2 + x_{2n - 1 - {i_j}}^2} \right)}  + l: =  - \varepsilon {x_N} + {C_5},
\end{equation}
where ${C_5}$ is finite. Hence
\[\frac{{x_{{i_1}}^2 + x_{{i_2}}^2}}{{ - {x_N}}} \to \varepsilon ,\]
and
\[\frac{{\left| {{x_{{i_1}}}} \right|}}{{ - {x_N}}} = \frac{{x_{{i_1}}^2}}{{ - {x_N}}}\frac{1}{{\left| {{x_{{i_1}}}} \right|}} \le \frac{{x_{{i_1}}^2 + x_{{i_2}}^2}}{{ - {x_N}}}\frac{1}{{\left| {{x_{{i_1}}}} \right|}} \to 0,\]
\[\frac{{\left| {{x_{{i_2}}}} \right|}}{{ - {x_N}}} = \frac{{x_{{i_2}}^2}}{{ - {x_N}}}\frac{1}{{\left| {{x_{{i_2}}}} \right|}} \le \frac{{x_{{i_1}}^2 + x_{{i_2}}^2}}{{ - {x_N}}}\frac{1}{{\left| {{x_{{i_2}}}} \right|}} \to 0.\]
It follows by $h(x) = 0$ that
\begin{equation}\label{eq37}
  {x_{{i_1}}}{x_{2n - 1 - {i_1}}} + {x_{{i_2}}}{x_{2n - 1 - {i_2}}} =  - \sum\limits_{j = 3}^l {{x_{{i_j}}}{x_{2n - 1 - {i_j}}}}  - \gamma x_N^3 - \frac{1}{2}\varepsilon {x_N}: =  - \gamma x_N^3 - \frac{1}{2}\varepsilon {x_N} + {C_6},
\end{equation}
where ${C_6}$ is finite, so
\[\frac{{{x_{{i_1}}}{x_{2n - 1 - {i_1}}}}}{{ - {x_N}}} + \frac{{{x_{{i_2}}}{x_{2n - 1 - {i_2}}}}}{{ - {x_N}}} =  - \gamma \frac{{x_N^3}}{{ - {x_N}}} - \frac{1}{2}\varepsilon \frac{{{x_N}}}{{ - {x_N}}} + \frac{{{C_6}}}{{ - {x_N}}} = \gamma \frac{{ - x_N^3}}{{ - {x_N}}} + \frac{1}{2}\varepsilon  + \frac{{{C_6}}}{{ - {x_N}}}.\]
The left hand side tends to 0 and the right hand side tends to $ - \infty $, a contradiction.

The remaining cases can be treated similarly. Then $\{ {x_N}\} $ must be bounded.

Proposition \ref{Pro21} is proved.

\begin{remark}\label{Re31}
a) In the proof of Proposition \ref{Pro21}, if ${x_{{i_1}}} \to \infty $ in 1) is changed to that one of ${x_{{i_2}}}, \cdots,{x_{{i_l}}},{x_{2n - 1 - {i_1}}}$, ${x_{2n - 1 - {i_2}}}, \cdots ,{x_{2n - 1 - {i_l}}}$ tends to $\infty$, then one can solve as in 1).

b) As a generalized case of 2) in the proof of Proposition \ref{Pro21}, if components
${x_i},{x_{2n - 1 - i}}(i \in {I_n})$ tend to $\infty$, then
\[\sum\limits_{i \in {I_n}} {\left( {x_i^2 + x_{2n - 1 - i}^2} \right)}  \to \infty .\]
It follows by $g(x) = 0$ that
\[\sum\limits_{i \in {I_n}} {\left( {x_i^2 + x_{2n - 1 - i}^2} \right)}  =  - \varepsilon {x_N} - l,\]
so
\[\frac{{\sum\limits_{i \in {I_n}} {\left( {x_i^2 + x_{2n - 1 - i}^2} \right)} }}{{ - {x_N}}} \to \varepsilon  > 0.\]
We have from $h(x) = 0$ that
\[\sum\limits_{i \in {I_n}} {{x_i}{x_{2n - 1 - i}}}  =  - \gamma x_N^3 - \frac{1}{2}\varepsilon {x_N} \to  - \infty ,\]
hence
\[\frac{{\sum\limits_{i \in {I_n}} {{x_i}{x_{2n - 1 - i}}} }}{{ - {x_N}}} = \frac{{ - \gamma x_N^3}}{{ - {x_N}}} - \frac{1}{2}\varepsilon \frac{{{x_N}}}{{ - {x_N}}} \to  - \infty ,\]
and by the Cauchy inequality,
\begin{align*}
   + \infty  \leftarrow \frac{{\left| {\sum\limits_{i \in {I_n}} {{x_i}{x_{2n - 1 - i}}} } \right|}}{{ - {x_N}}}& \le \frac{{{{\left( {\sum\limits_{i \in {I_n}} {x_i^2} } \right)}^{\frac{1}{2}}}{{\left( {\sum\limits_{i \in {I_n}} {x_{2n - 1 - i}^2} } \right)}^{\frac{1}{2}}}}}{{ - {x_N}}} \le \frac{{{{\left( {\sum\limits_{i \in {I_n}} {\left( {x_i^2 + x_{2n - 1 - i}^2} \right)} } \right)}^{\frac{1}{2}}}{{\left( {\sum\limits_{i \in {I_n}} {\left( {x_i^2 + x_{2n - 1 - i}^2} \right)} } \right)}^{\frac{1}{2}}}}}{{ - {x_N}}} \\
   &  = \frac{{\sum\limits_{i \in I} {\left( {x_i^2 + x_{2n - 1 - i}^2} \right)} }}{{ - {x_N}}} \to \varepsilon ,
\end{align*}
a contradiction.

c) To the generalized case of 3) in the proof of Proposition \ref{Pro21}, if ${x_i}(i \in {I_n})$ tend to $\infty$ and ${x_{2n - 1 - i}}(i \in {I_n})$ are bounded, then
\begin{center}
$\sum\limits_{i \in {I_n}} {x_i^2}  \to \infty $ and $\sum\limits_{i \in {I_n}} {x_{2n - 1 - i}^2} $ is bounded.
\end{center}
It uses $g(x) = 0$ to have
\[\sum\limits_{i \in {I_n}} {x_i^2}  =  - \sum\limits_{i \in {I_n}} {x_{2n - 1 - i}^2}  - \varepsilon {x_N} + l,\]
and
\[\frac{{\sum\limits_{i \in {I_n}} {x_i^2} }}{{ - {x_N}}} \to \varepsilon  > 0.\]
It follows from $h(x) = 0$ that
\[\sum\limits_{i \in {I_n}} {{x_i}{x_{2n - 1 - i}}}  =  - \gamma x_N^3 - \frac{1}{2}\varepsilon {x_N} \to  - \infty ,\]
and
\[\frac{{\sum\limits_{i \in {I_n}} {{x_i}{x_{2n - 1 - i}}} }}{{{{\left( { - {x_N}} \right)}^{\frac{1}{2}}}}} = \frac{{ - \gamma x_N^3}}{{{{\left( { - {x_N}} \right)}^{\frac{1}{2}}}}} - \frac{1}{2}\varepsilon \frac{{{x_N}}}{{{{\left( { - {x_N}} \right)}^{\frac{1}{2}}}}} \to  - \infty ,\]
then by the Cauchy inequality,
\[ + \infty  \leftarrow \frac{{\left| {\sum\limits_{i \in {I_n}} {{x_i}{x_{2n - 1 - i}}} } \right|}}{{{{\left( { - {x_N}} \right)}^{\frac{1}{2}}}}} \le \frac{{{{\left( {\sum\limits_{i \in {I_n}} {x_i^2} } \right)}^{\frac{1}{2}}}{{\left( {\sum\limits_{i \in {I_n}} {x_{2n - 1 - i}^2} } \right)}^{\frac{1}{2}}}}}{{{{\left( { - {x_N}} \right)}^{\frac{1}{2}}}}} \to \sqrt \varepsilon  {\left( {\sum\limits_{i \in {I_n}} {x_{2n - 1 - i}^2} } \right)^{\frac{1}{2}}},\]
a contradiction.
\end{remark}

\textbf{Proof of Proposition \ref{Pro22}} Let us apply the contradiction. Assume that the rank of the Jacobian matrix for $g(x)$ and $h(x)$ is smaller than 2, then there exists $\kappa  \ne 0$, such that
\begin{equation}\label{eq38}
\nabla g(x) = \kappa\nabla h(x).
\end{equation}
For $i \in I_n$, it has ${g_{{x_i}}} = 2{x_i},\;{g_{{x_{2n - 1 - i}}}} = 2{x_{2n - 1 - i}},\;{h_{{x_i}}} = {x_{2n - 1 - i}},\;{h_{{x_{2n - 1 - i}}}} = {x_i},$ and by \eqref{eq38} that
\begin{equation}\label{eq39}
\left\{ {\begin{array}{*{20}{c}}
{2{x_i} = \kappa {x_{2n - 1 - i}},}\\
{2{x_{2n - 1 - i}} = \kappa {x_i},}
\end{array}} \right.
\end{equation}
i.e.,
$\left( {1 - {{\left( {\frac{\kappa }{2}} \right)}^2}} \right){x_i} = 0,$
then
\begin{center}
${x_i} = 0$ or $\kappa  = 2, - 2.$
\end{center}

For $i = N$, it follows ${g_{{x_N}}} = \varepsilon ,$ ${h_{{x_N}}} = 3\gamma x_N^2 + \frac{1}{2}\varepsilon ,$ and by \eqref{eq38} that
\begin{equation}\label{eq310}
3\kappa \gamma x_N^2 =  - \frac{1}{2}\kappa \varepsilon  + \varepsilon .
\end{equation}

We can show that all cases above yield contradictions. Actually, when ${x_i} = 0$, we have ${x_{2n - 1 - i}} = 0$ from $2{x_{2n - 1 - i}} = \kappa {x_i}$ in \eqref{eq39}, and so
\[0 = g(x) = \varepsilon {x_N} - l,\]
\[0 = h(x) = \gamma x_N^3 + \frac{1}{2}\varepsilon {x_N}.\]
It yields a contradiction as in Case ${a_1}),{b_1})$.

When $\kappa  = 2$, we have from \eqref{eq39} and \eqref{eq310} respectively that ${x_i} = {x_{2n - 1 - i}}$ and $x_N^2 = \frac{{ - \frac{1}{2}\kappa \varepsilon  + \varepsilon }}{{3\kappa \gamma }} = \frac{{ - \varepsilon  + \varepsilon }}{{6\gamma }} = 0,$ i.e., ${x_N} = 0$, then
\[0 = h(x) = \sum\limits_{i \in {I_n}} {x_i^2} ,\]
and
$${x_i} = {x_{2n - 1 - i}} = 0,$$
so
$$0 = g(x) =  - l < 0,$$
a contradiction.

When $\kappa  = -2$, it yields from \eqref{eq39} and \eqref{eq310} respectively that ${x_i} =  - {x_{2n - 1 - i}}$ and
$$x_N^2 = \frac{{\varepsilon  + \varepsilon }}{{ - 6\gamma }} =  - \frac{\varepsilon }{{3\gamma }} = {l^{\frac{2}{3}}},$$
then
\begin{equation}\label{eq311}
{x_N} =  \pm {l^{\frac{1}{3}}}.
\end{equation}
Using
\[0 = g(x) = 2\sum\limits_{i \in {I_n}} {x_i^2}  + \varepsilon {x_N} - l,\]
\[0 = 2h(x) =  - 2\sum\limits_{i \in {I_n}} {x_i^2}  + 2\gamma x_N^3 + \varepsilon {x_N},\]
we have
\[2\gamma x_N^3 + 2\varepsilon {x_N} - l = 0,\]
and obtain as in Case ${a_3}),{b_1})$ that
\[{x_N} =  \left( { - {{\left( {\sqrt 2  - 1} \right)}^{\frac{1}{3}}} - {{\left( {\sqrt 2  + 1} \right)}^{\frac{1}{3}}}} \right){l^{\frac{1}{3}}}.\]
It is different from ${x_N}$ in \eqref{eq311}), a contradiction.

Proposition \ref{Pro22} is proved.

\section{Proof of Theorem \ref{Th1} (even integers) }

For the even integers, supposing that the Cai-Zhang-Shen conjecture is not true, then there exists an even integer $2n$ such that $2n$ can not be expressed as the sum of two figurate primes. Let us take respectively
\[P = \left( {\delta ({i_1}), \cdots ,\delta ({i_l}),\delta (n),\delta (2n - {i_l}), \cdots ,\delta (2n - {i_1})} \right),\]
$$f(x) = \sum\limits_{i \in {I_n}} {{x_i}{x_{2n - i}}}  + s{x_n^2},$$
\[g(x) = \sum\limits_{i \in {I_n}} {\left( {x_i^2 + x_{2n - i}^2} \right)}  + \varepsilon {x_n} - l,\]
\[h(x) = \sum\limits_{i \in {I_n}} {{x_i}{x_{2n - i}}}  + \gamma x_n^3 + \frac{1}{2}\varepsilon {x_n}.\]
Similarly to the proof for odd integers in Sections 2, we also reach a contradiction.

\section{Conclusions}
In previous sections, we prove Cai-Zhang-Shen conjecture for figurate primes.
The way of proof  really provides a new approach to confirm Goldbach's binary
conjecture. It is worth trying and we will further consider the well known and difficult
conjecture.

\textbf{Acknowledgments.} We are especially indebted to the anonymous referees for the careful readings and many useful suggestions.

\textbf{Conflicts of Interest} The authors declare that there is no conflict of interest regarding the publication of this paper.

\end{document}